\documentstyle[12pt]{amsart}

\newtheorem{thm}{Theorem}[section]

\newtheorem{exa}[thm]{Example}
\newtheorem{propn}[thm]{Proposition}
\newtheorem{rem}[thm]{Remark}
\newtheorem{lemma}[thm]{Lemma}
\newtheorem{cor}[thm]{Corollary}

\def\subsetneqq{\subset}

\def \mo{{\cal M}(n,d)} 
\def \bn{{\cal W}^{k-1}_{n,d}} 
\def \bnp{{\cal W}^{k-1}_{n,d''}} 
\def \bns{{\cal W}^{sk-1}_{n,d}} 
\def \bnl{{\cal W}^{k-1}_{1,d}} 
\def \bnn{{\rho}^{k-1}_{n,d}} 
\def \pa{(\mu,\lambda)} 
\def \bnnl{{\rho }^{k-1}_{1,d}} 
\def \wmo{\widetilde{\cal M}(n,d)} 
\def \wbn{\widetilde{\cal W}^{k-1}_{n,d}} 
\def \parti{\eta(s)} 
\def \hparti{\hat\parti} 
\begin{document}
\baselineskip=15.5pt

\title[Nonemptiness of Brill-Noether]{Nonemptiness of 
Brill-Noether loci} 

\author{L. Brambila-Paz}

\address{CIMAT, Apdo. Postal 402,
C.P. 36240. Guanajuato, M\'exico}

\email{lebp@@fractal.cimat.mx}

\author{ V. Mercat}

\address{Institut de Mathematiques
Universite Paris 7-case 7012, 
Jussieu 75251,  Paris cedex 05, France}

\email{mercat@@math.jussieu.fr}

\author{P. E. Newstead}

\address{Department of Mathematical Science, The University of
 Liverpool, Liverpool, L69 3BX, England}

\email{newstead@@liverpool.ac.uk}

\author{F. Ongay}

\address{CIMAT, Apdo. Postal 402,
C.P. 36240. Guanajuato, M\'exico}
 
\email{ongay@@fractal.cimat.mx}

\begin{abstract}Let $X$ be a non-singular algebraic
 curve of genus $g$. 
We prove that the Brill-Noether locus $\bns $ 
is non-empty
if $d= nd' +d'' $ with $0< d'' <2n$, $1\le s\le g$,
$d'\geq (s-1)(s+g)/s $, $n\leq d''+(n-k)g$,
$(d'',k)\ne(n,n)$. These results hold for an 
arbitrary curve of genus $\ge 2$,
and allow us to construct  a region in the associated
``Brill-Noether $\pa$-map'' of points for which the
Brill-Noether loci are non-empty. Even for the generic
case, the region so constructed extends beyond that
defined by the so-called ``Teixidor parallelograms.''
For hyperelliptic curves, the same methods give more
extensive and precise results.
\end{abstract}

\thanks{
authors are members of the Europroj research group Vector
Bundles on Algebraic Curves. The first and fourth authors
acknowledge support from CONACYT (grant nos. 28-491-E and 28-492-E) } 

\date{}

\maketitle

\section{Introduction}

Brill-Noether theory is concerned with the 
study of the subvarieties of
the moduli space of stable bundles, determined 
by bundles having at least a
specified number of sections. More precisely, if
 $\mo$ is the moduli space of stable
vector bundles  of rank $n$ and degree $d$ over 
 a non-singular algebraic curve $X$ of genus
$g\geq 2$ over {\bf{C}},  and $k\geq 1$,  the corresponding {\it
Brill-Noether locus} is
$$\bn := \{ E\in \mo | h^0(E)\geq k
\}. $$
 The main questions in Brill-Noether theory regard 
the nonemptiness, dimension, connectedness,
irreducibility, cohomology classes, etc., of these varieties.
(Similar statements can be made for semistable bundles.)

For line bundles,  Brill-Noether theory has been 
studied since the last
century, and for a generic curve the basic questions 
have been answered (see \cite{ACGH}). However, the 
corresponding theory for vector bundles of higher rank
is far from being complete even for the generic case. 
(In section 2 we will recall the known results for this case.)

In this paper we will be concerned with the 
nonemptiness question.
 We will prove (see Corollary 3.2 and Theorem 3.9) that 

{\it if $d=nd'+d''$, $0<d''<2n$,
 and $d' \ge0 $ then, {\it for any}  $X$, $$\bn  \neq \emptyset \
\ \ \ \  \text{if} \ \ \ \  n\leq d''+(n-k)g\text{ and }
(d'',k)\ne(n,n).\eqno(A)$$}
More generally, 

{\it if $1\le s\le g$ and there exists a
line bundle $L$ on $X$ of degree $d'$ with $h^0(L)\ge
s$, then    $${\cal W}^{sk-1}_{n,d} \neq \emptyset \ \ \
\text{if}  \ \ \  0<d''<2n,\ n\leq d''+(n-k)g \text{ and }
(d'',k)\ne(n,n).\eqno(B)$$ }

Observe that such a line bundle $L$ always exists 
if $d'\geq
(s-1)(s+g)/s$. If $X$ is hyperelliptic, we see by
considering powers of the hyperelliptic line bundle 
that $L$ exists if $d' \geq 2s-2.$

To get a clear picture of the triples $(n,d,k)$ 
for which $\bn$
is not empty, as Alastair King has pointed out, 
it is easier to 
represent
this kind of result if we write $\mu=d/n$,
$\lambda=k/n$ and plot points in the $\pa$-plane; 
we will refer
to this representation as the {\it Brill-Noether map} 
(or BN
map). It
will be convenient to call a point $\pa\in{\bold Q}^2$
a {\it $n$-Brill-Noether} point (or $n$-BN point),
 $n\in\Bbb N$,
if $d=n\mu$ and $k=n\lambda$ are both integers and 
$\bn\ne\emptyset$. If $\pa$ 
is $n$-BN for all $n$ such that  $d=n\mu$ and 
$k=n\lambda$ are both integers, we just say it is a BN point.

The results $(A)$ and $(B)$ define a
region in the $\pa$-plane, which we denote by $BMNO$,  
 where the points are
$n$-BN for ``many values'' of $n$, and thus the 
corresponding
Brill-Noether loci are non-empty (for an explanation 
of what we mean
by ``many values'' see Remarks 4.3 and 4.4). Actually, in the
 hyperelliptic case we give a precise description 
of which points are $n$-BN.
 The results, in this case, 
come close to a complete solution of the
nonemptiness problem. In particular, the boundary of the
region in which stable bundles of rank $>1$ can
exist is completely determined, and, as might be
expected, it is close to the Clifford line.

 On the other hand, the results in \cite{Te1}, 
later refined in \cite{M2},
show that one can define a polygonal region $T$, the so-called
``Teixidor's parallelograms,"  such  that all the 
points $\pa$ inside $T$
are BN, except perhaps at certain vertices. The region
$BMNO$ covers a
large part of $T$, but more importantly it extends beyond
$T$. Our methods, and especially Theorem 3.9, give stronger
results for special curves. Furthermore, since we do 
not use the results of \cite{Te1} and \cite{M2}, 
our results give another proof of
nonemptiness for those parts of $T$ which are included in
$BMNO$.

 \vskip 6pt 
The paper is organized as follows:
In section 2 we give a brief  survey of what is known about
$\bn$. In section 3,  we prove assertions $(A)$ and $(B)$.
In section 4 we describe the region $BMNO$ in the
$\pa$-plane for the general case. In section 5 
we compare the regions
$BMNO$ and $T$. In section 6 we study the case where $X$
is hyperelliptic.

\section{A survey of the known theory}

In this section we recall the known results of 
Brill-Noether theory for
vector bundles of higher rank (see also \cite{BGN} and
\cite{M1}), and thereby also fix
 notations. We will use the Brill-Noether map to 
indicate the  regions 
where not only
 the nonemptiness is known but also some of the topology of $\bn$.

\bigskip

Let $X$ be a non-singular algebraic curve of genus $g\geq
2$ over {\bf{C}} and $\mo$ the moduli space of stable
vector bundles over $X$ of rank $n$ and degree $d$. We
define the Brill-Noether loci for $(n,d,k)$, 
$$\bn := \{ E\in \mo |
h^0(E)\geq k \} $$  as in the introduction. Denoting by
$\wmo$ the moduli space of equivalence classes of
semistable bundles over $X$ of rank $n$ and degree $d$,
we can define similarly
 Brill-Noether loci $\wbn$ in $\wmo$. In what follows we
will concentrate on stable bundles, but we will also
explicitly indicate where semistable bundles are allowed.
Since for $k>0,d<0$, $\bn =\emptyset$ and for $ k\leq 0, \  \bn $ is 
the whole moduli space $\mo $, we will 
assume that $d \geq 0$ and $k \geq 1.$

 It follows from the theory of determinantal varieties
that every non-empty component of $\bn$ has dimension greater
than or equal to the Brill-Noether number
$$\bnn := n^2(g-1)+1-k(k-d+n(g-1)),$$ and for generic $X$
this number is the expected  dimension of $\bn$. However, 
there is
no similar formula for the expected dimension of $\wbn$.

For
$n=1$, it is a classical result (see \cite{ACGH}) that 
$\bnl\ne\emptyset$ if $ \bnnl\geq 0$; moreover, for a
generic curve the converse is also true. However,
for $n\geq 2$, it is known (see \cite{BGN}) that $\bnn
\geq 0$ does not imply that $\bn\ne\emptyset$.

A point $\pa\in{\bold Q}^2$ will be called
a {\it $n$-Brill-Noether} point (or $n$-BN point),
 $n\in\Bbb N$,
if $d=n\mu$ and $k=n\lambda$ are both integers and 
$\bn\ne\emptyset$. If
it is $n$-BN for all $n$ such that $d=n\mu$ and 
$k=n\lambda$ are both integers we just say it is a BN point.

 In the Brill-Noether map,
using the Riemann-Roch Theorem and
Clifford's Theorem, one can define a region such that
outside this region the problem becomes trivial, in the
sense that $\bn$ is either empty or the whole moduli space
$\mo$. More precisely,   
consider the following lines:

\vskip 4pt

{\it i\/}) $\mu =\lambda + g-1$ \ \ \  (Riemann-Roch line)

\vskip 4pt

{\it ii\/}) $\mu = 2\lambda -2 $  \ \ \ \ (Clifford line )

\vskip 4pt

These lines, together with the positive axes and the line
$\mu=2g-2$, define a bounded pentagonal region, which we
denote by $P$. That is, we define $P$ to be the
region defined by the inequalities
$$\mu<\lambda+g-1,\ \mu\geq
2\lambda-2,\ 0\leq\mu\leq2g-2,\ \lambda>0$$
(see Figure 1).

 Below and to the
right of $P$, Riemann-Roch implies that $\bn$ is the whole
space. 
Above and to the left of $P$, Riemann-Roch and
Clifford's  Theorem, together with the definition of
stability, imply that $\bn$ is empty. 
  Thus, we are interested in studying only the points
inside $P$.

\begin{rem}\begin{em} {\it i)} For
$\mu=0$, the only stable bundle in $P$  is the trivial
line bundle $\cal O$ at the point $(0,1)$, while for
$\mu=2g-2$, the only such bundle is the canonical line
bundle $K$ at the point
$(2g-2,g)$. However semistable bundles exist at all
points of these two edges of $P$ \cite{BGN}.
 Note that, according to our definition, the points
$(0,1)$ and $(2g-2,g)$ are only $1$-BN.

\noindent {\it ii)} The
inequalities defining $P$ are all sharp except for the
Clifford bound. The exact r\^ole of this bound is not
clear but it has, for example, been improved
 by Re (see \cite{R}) for non-hyperelliptic curves. In
this case, if we restrict to the range
$1\leq\mu\leq2g-3$, we can replace the Clifford line by
the line $\mu =2\lambda -1$. For further improvements,
see for example \cite{BN} and \cite{M4}.
\end{em}\end{rem}

By Serre duality we know that, if $\pa$ is BN, 
then so is $$\sigma\pa:=(2g-2-\mu,\lambda +g-1-\mu).$$
Though it is not readily apparent from Figure 1, this
gives a symmetry in $P$ through
  the line $\mu = g-1$. For later purposes, it will 
be convenient to write 
$$R=\{\pa\in P: \mu\le g-1\},$$ 
so that in particular $P=R\cup\sigma(R)$ (see Figure 1).

An important feature of the BN map is the curve
defined by the equation
$$\widetilde{\rho}=\frac{1}{n^2}(\bnn -1) =0,$$ called the
 {\it Brill-Noether curve} (or BN curve). From what was
said earlier, this represents
 the boundary of the region where one would expect the
Brill-Noether loci to have positive dimension, though it is
 known that this analogy to the case of line bundles is not
valid in
 general (see \cite{BGN}, \cite{M1}, etc.). The BN
curve is a portion of a hyperbola, with equation
$$\widetilde{\rho}\pa =(g-1)-\lambda (\lambda -\mu
+g-1)=0.$$

The results of \cite{Te1} and \cite{M2} allow us to define a
polygonal region which we denote by  $T$, contained in
 the interesting region  $P$ in the $\pa$-plane, such that
any point $\pa$ in $T$ is BN except  possibly for certain
vertices. This region was described in detail in
\cite{BGN} and \cite{M1} in its original form, and we will
recall its construction in section 4 incorporating the
results of \cite{M2}. For the time being, we just point
out that $T$ has  sides parallel  to the lines $\lambda
=0,\mu =\lambda$, and vertices at points with integer
 coordinates, on or below the BN curve (see
Figure 2).

 The most significant results for our purposes
are the following, which hold for slopes restricted to
$0\leq\mu < 2$:

\begin{enumerate}

\item{}  For $0< \mu \leq 1$, Brambila-Paz,
Grzegorczyk and Newstead proved in \cite{BGN} that
$\pa$ is BN if and only if $1\leq\mu+ (1-\lambda )g $ and
$(\mu,\lambda)\ne(1,1)$.

\item{}   For $1< \mu <2$, Mercat in
\cite{M1} proved that $\pa$ is BN if and only
if $1\leq \mu+ (1-\lambda )g $. 

\end{enumerate}

In the BN map these
results define two trapezoidal regions inside $R$,
 which we will denote by $BGN$ and $M$ respectively (see
Figure 3).

Ballico, Mercat and Newstead have
recently proved the existence of stable bundles at some
points outside the regions defined above; in particular,
these bundles can have negative Brill-Noether number
\cite{BMN}.

 For $X$ generic, Teixidor also proved that $\bn$
has an irreducible component of dimension $\bnn$.
For any curve, Brambila-Paz,
Grzegorczyk and Newstead in \cite{BGN} also proved that
 for $ 0<\mu \leq 1$, if
$\bn$ is nonempty, then it is irreducible of dimension
$\bnn$, and $\text{Sing\,}\bn ={\cal W}^k_{n,d}$
(\cite{BGN}, Theorem A). For $1<\mu <2$, Mercat in
 \cite{M1} also proved that, if $n=d+(n-k)g $ or
$n<d<n+g$, $\bn $ is irreducible (\cite{M1}, 2-B-1 and
3-A-1); in any case all components have the expected
dimension and $\text{Sing\,}\bn ={\cal W}^k_{n,d}$
(\cite{M1}, 2-C-1).
 So, for slopes $0\leq \mu < 2$, the
results are very complete.

For $k=1$, Sundaram \cite{Su} proved that ${\cal
W}^0_{n,d}$ is irreducible of dimension $\rho^0_{n,d}$,
and Laumon \cite{L} showed $\text{Sing\,}{\cal
W}^0_{n,d}$= ${\cal W}^1_{n,d}$. For rank $2$ there are 
also results of nonemptiness and irreducibility in \cite{T}, 
\cite{Te2} and \cite{Su} and for rank $3$ in \cite{NB}.

\begin{rem}\begin{em} If $X$ is not
hyperelliptic, the results of \cite{M1} can be extended
to cover the case $\mu=2$ \cite{M3}. For further
details, see Remark 4.8.
\end{em}\end{rem}

>From 
the results of \cite{BGN} and  the symmetry of the region
$P$, for $g=2$ one has a complete description  of
Brill-Noether loci: nonemptiness, irreducibility and
singularities.

 \section{ Nonemptiness of Brill-Noether loci}

In this section we will prove assertions $(A)$ and $(B)$.

\smallskip

Note first that, for any line bundle $L$ of degree
$d'$, the formula  $E\mapsto E\otimes L $ defines an
isomorphism $$ \Phi _L :{\cal M}(n,d) \rightarrow  {\cal
M}(n,nd'+d). \eqno(1)$$

We shall make repeated use of this idea of
tensoring stable bundles by line bundles. If $V$ is a
subvariety of $\mo$, then $\Phi _L(V)$ is a subvariety
of ${\cal
M}(n,nd'+d)$; in particular, if $V$ is a non-empty
Brill-Noether locus and $h^0(L)>0$, then $\Phi _L(V)$ will 
meet certain
Brill-Noether loci in ${\cal M}(n,nd'+d)$, which will
therefore be non-empty.

Actually, if
$d'\ge0 $ then there always exists
 a line bundle $L$ of degree $d'$ that has at least one
section, so we have the
following theorem.

\begin{thm} If $\bn \neq \emptyset$
(respectively
  $\wbn\neq\emptyset$), then ${\cal W}^{k-1}_{n,d+nd'}\neq
  \emptyset$ (respectively $\widetilde{\cal
W}^{k-1}_{n,d+nd'}\neq\emptyset$) for any $d'\ge0$.
\end{thm}

\smallskip

{\it Proof:} Choose a line bundle $L$ of degree $d'$ with 
$h^0(L) >0$. Then, for any $E$ with $h^0(E)\ge k$ 
 we have $h^0(E\otimes L)\ge k$.

\bigskip

\begin{cor}{\bf (Assertion (A))} Suppose $d=nd'+d''$ with
$0<d''<2n$, $d'\ge0$, $n\leq d''+(n-k)g$ and
$(d'',k)\ne(n,n)$. Then $\bn\ne\emptyset$.
\end{cor}

\smallskip

{\it Proof:} This follows from the results of \cite{BGN}
and \cite{M1} and Theorem 3.1. 

\bigskip

\begin{cor} If $k<n$, then ${\cal
W}^{k-1}_{n,rn} \neq \emptyset
  $ for $r\geq 1$.
\end{cor}

\smallskip

{\it Proof:} Take $d''=n$ and $d'=r-1$ in Corollary 3.2.

\bigskip

 In general, the multiplication map  $$ \mu _{E,F}
: H^0(E)\otimes H^0(F) \rightarrow H^0(E\otimes F) $$ is
not injective and $h^0(E)^.h^0(F)$ does not give a lower
bound for $h^0(E\otimes F)$. However, if $E$ is a 
point in $\bn$ with $d/n <2$ and $L$ a line bundle of
degree $d'\geq0 $ with at least $s$ independent 
sections, then we have
the following lemmas.

\bigskip

\begin{lemma} If $d<n+g$, then $h^0(E\otimes L) \geq
ks.$ 
\end{lemma}

\smallskip

{\it Proof:} From \cite{BGN} (for $d\le n$) and \cite{M1},
3-A-1 (for $n<d<n+g$) we know that any such bundle has
${\cal O}^k $ as a subsheaf. Hence, $\oplus ^kL $ is a
subsheaf of $E\otimes L$. Therefore  $h^0(E\otimes L) \geq
h^0(\oplus ^kL) \geq ks.$

\bigskip

\begin{rem}\begin{em} Note that the
existence of $E$ implies that $n\le d+(n-k)g$, so the
hypothesis $d<n+g$ implies that $k\le n$.
\end{em}\end{rem}

\bigskip

\begin{lemma} If $k>n$, $d=n +(k-n)g$ and $d'\le2g$, then 
$h^0(E\otimes L) \geq ks.$
\end{lemma}

\smallskip

{\it Proof:} From  \cite{M1}, 2-B-1 we know that any such
bundle fits in an exact sequence  $$0\rightarrow F^*
\rightarrow {\cal O} \otimes H^0(X,F)^* \rightarrow
E\rightarrow 0\eqno(2)$$ where $F$ is a stable bundle of
slope $>2g$ and $h^0(F)=h^0(E)=k$. Tensor
$(2)$ by $L$ and take the cohomology sequence. Since
$\deg(F^*\otimes L)<0$, we have  $H^0(X,F^*\otimes
L) =0$ and hence
 $$h^0(E\otimes
L) \geq h^0(L)^.h^0(F) \geq ks.$$ 

\bigskip

\begin{lemma} If $n+g\le d<2n$ and $d'\le2g$, then
there exists $E\in\bn$ with $h^0(E\otimes L) \geq ks.$
\end{lemma}

\smallskip

{\it Proof:} In this case we have two exact sequences
$$0\rightarrow {\cal O}^{l'} \rightarrow E'\rightarrow
E\rightarrow 0\eqno(3)$$ $$0\rightarrow D(E')^*
\rightarrow {\cal O}^{n+l+l'} \rightarrow E' \rightarrow
0\eqno(4)$$
where $k\le n+l$, $E\in{\cal W}_{n,d}^{n+l-1}$, $E'$ and
$D(E')$ are stable and $\mu(D(E'))>2g$ (see \cite{M1}, 
3-B-1 and
its proof).
 Tensor both sequences by $L$ and take the cohomology
sequences. Since $H^0(X,L\otimes D(E')^*)=0$,
$$h^0(E'\otimes L) \geq h^0(\oplus ^{n+l+l'}L)=
(n+l+l')h^0(L)$$ Thus $$h^0(E\otimes L)\ge(n+l)h^0(L) \geq
ks.$$

\bigskip
\begin{rem}\begin{em}{\it i)} From the
Brill-Noether theory for line bundles we know that, if
$d'\geq\parti:= (s-1)(s+g)/s$, then there exists a line
bundle $L$ of degree $d'$ with at least $s$  independent 
sections.
Actually, ${\cal W}^{s-1}_{1,d'}$ is the variety of such
bundles and, for $X$ generic, it has dimension
$g-s(s-d'+g-1)=s(d'-\eta(s))$.

\noindent {\it ii)} 
If $X$ is hyperelliptic, we see by
considering powers of the hyperelliptic line bundle 
that $L$ exists if $d' \geq 2s-2.$ 
\end{em}\end{rem}

  We deduce the following theorem that proves
assertion ($B$):

\bigskip

\begin{thm} Suppose $d=nd'+d''$ with $d'\ge0$,
  $0<d''<2n$ and that $1\le s\le
g$. If $n\leq d'' +(n-k)g$ and
$(d'',k)\ne(n,n)$ and there exists a line bundle 
$L$ on $X$ of degree
$d'$ with $h^0(L)\ge s$, then  $\bns \neq \emptyset$.
\end{thm}

\smallskip

{\it Proof:} Note first that, for fixed $n$, $k$,
$d''$, $s$, if the theorem is true for one value of
$d'$, it is true for all larger values. Since $\eta$
is an increasing function of $s$ and $\eta(g)=2g-2$,
it is therefore sufficient by Remark 3.8 to prove the
theorem with the additional hypothesis that $d'\le2g-2$.
Now, from \cite{BGN} and \cite{M1} we know that under the
given hypotheses $\bnp$ is non-empty. 
It follows from the lemmas that $\bns $ is
non-empty. 

\bigskip

\begin{rem}\begin{em} In the
semistable case, the theorem can be extended to the cases
$d''=0$, $k\le n$ and  $(d'',k)=(n,n)$.
\end{em}\end{rem}

\bigskip

\begin{cor} If $k<n$ and $r\geq s+g-g/s $, then ${\cal
    W}^{sk-1}_{n,rn}\neq\emptyset.$ 
\end{cor}

\smallskip

{\it Proof:} Take $d''=n$ and $d'=r-1$ in Theorem 3.9.

\bigskip

We finish this section by describing the above results for $g=3$.

\bigskip

\begin{exa}\begin{em} Suppose $X$ has
genus $3$. It follows from \cite{BGN}, \cite{M1} and
Corollary 3.3 that a point $\pa\in R$ is BN if
$\mu>0$, $1\le \mu+3(1-\lambda)$, $\pa\ne(1,1)$, except
possibly when $\mu=2$ and $1\le \lambda\le 4/3$. In fact
$(2,1)$ is also BN by \cite{Te1} and \cite{M2}. Moreover, any
$\bn$ corresponding to such a point has pure dimension
$2n^2+1-k(k-d+2n)$ and    $\text{Sing\,}\bn ={\cal
W}^k_{n,d}$; it is irreducible if  $0<d<n+3$ or if
$d=n+3(k-n)$. 

If $X$ is not hyperelliptic, then by \cite{M3} we can
remove the exceptional case and say that $\pa\in R$ is BN
 if and only if
$\mu>0$, $1\le \mu+3(1-\lambda)$, $\pa\ne(1,1)$. The
statements about dimension and singularities still apply.
Apart from the trivial bundle $\cal O$, there is also
precisely one further stable bundle on $X$ in $R$, namely
the bundle $E_K$ with $n=2, d=4,
k=3$ (see \cite{M1}, 2-A-4).
\end{em}\end{exa}

\section{ The region $BMNO$ in the
$(\mu,\lambda$)-plane}

In this section,  using the
results of the previous section, combined with 
Serre duality, we describe the region $BMNO$. 
Throughout the section $X$ is an arbitrary 
non-singular algebraic curve of genus $g\geq 3$.

\bigskip

In order to translate the results of section 2 into
geometric form, we introduce for each $d'$ and $s$ such that 
$d'\ge\parti=(s-1)(s+g)/s $, the
affine map $T_{d',s}$ given by
$$T_{d',s}(\mu,\lambda)=(\mu+d',s\lambda).$$  
Notice that it shifts points to the right, and, if $s>1$,
also expands in the $\lambda$ direction. (We shall refer
to these maps as {\it translations} although strictly
speaking only the
$T_{d',1}$ are translations.) The  idea is to use the
regions $BGN$ and $M$  as
 ``tiles'' to cover a larger region, the tiling
being  obtained by translating $BGN$ and $M$ by the maps
$T_{d',s}$. Then we will apply Serre duality to obtain the
 $BMNO$ region. 

More precisely, recall that  $\parti=(s-1)(s+g)/s$, and set
$\hparti=\lceil\eta(s)\rceil$ (we will use the notation
$\lceil\cdot\rceil$, and $\lfloor\cdot\rfloor$,
respectively, for the least integer not smaller, and the
largest integer not greater than a given number, the
so-called ``ceiling'' and ``floor'' functions). We will
 consider $d'$ and $s$ such that $d' \geq \hparti$ for 
the affine maps $T_{d',s}$.  For the description of 
the regions 
 $BGN$ and $M$, we consider the
trapezia $$ BGN'=\{\pa: 0<\mu\leq1,\quad
0<\lambda\leq\frac{1}{g}(\mu+g-1)\}$$
$$M'=\{\pa:
1<\mu<2,\quad 0<\lambda\leq\frac{1}{g}(\mu+g-1)\}.$$ 
Note however that the point $(1,1)\in
BGN'$ is only $1$-BN, so we define $$BGN=BGN'-\{(1,1)\}.$$
Moreover, a geometrical interpretation of Corollary 3.3
shows that, if we translate $BGN$ by $T_{1,1}$, we
obtain the points in the boundary of  $M'$ with $\mu=2,
0<\lambda<1$, and so these also give BN points. We
therefore define
$$M=M'\cup\{(2,\lambda): 0<\lambda<1\}.$$

\begin{rem}\begin{em} Notice that
there are still  some points in the boundary line $\mu=2$
of $M$ that are not covered this way, namely those for 
which
$1\leq\lambda\leq 1+1/g$, and therefore we do not know
whether they are BN points or not. For the non-hyperelliptic case see
Remark 4.3. 
\end{em}\end{rem}

>From Theorem 3.9 and Remark 3.8 we have 

\begin{thm} If $\pa\in BGN\cup M$, $1\le s\le
g$ and $d'\ge\hparti$, then $T_{d',s}\pa$ is $n$-BN for all $n$ 
such that $(\mu ,\lambda )$ is $n$-BN.
\end{thm}

\begin{rem}\begin{em} One can in fact obtain many
points in $T_{d',s}(BGN\cup M)$  which are BN points. 
 Let
$(\mu,\lambda)\in T_{d',s}(BGN\cup M)$ and write
$\mu=\frac{a}{b}$, $\lambda=\frac{c}{e}$ in their lowest
terms. To get bundles of rank $n$, the conditions we need
are that $b|n$ and $es|cn$. If $s|c$, these conditions
reduce to $b|n$ and $e|n$, so $\pa$ is BN. Points of this
form are dense in $T_{d',s}(BGN\cup M)$.
\end{em}\end{rem}

\begin{rem}\begin{em} {\it i)}
 If we take a point $\pa$ in
$T_{d',s}(BGN\cup M)$ lying strictly below the top
boundary, we can obtain an improvement to Theorem 4.2.
 Suppose $n$ is a
positive integer such that $n\mu$ and $n\lambda$ are both
integers and define
$\lambda'=\frac{s}{n}\lceil\frac{n\lambda}{s}\rceil$. If
$(\mu,\lambda')\in T_{d',s}(BGN\cup M)$, then ${\cal
W}^{n\lambda'-1}_{n,n\mu}\ne\emptyset$; hence ${\cal
W}^{n\lambda-1}_{n,n\mu}\ne\emptyset$. Now
$$\lambda'\le\frac{s}{n}\left(\frac{n\lambda+s-1}{s}
\right)=\lambda+\frac{s-1}{n};$$
so this holds for all sufficiently large $n$ (even if 
$n\lambda/s$ is not an integer).

\noindent {\it ii)} 
By a similar method, taking $\lambda ' =s$, we can 
show also that the region
$$\mu>\hparti+1,\ \lambda\le s,\ \mu\not\in{\bold N}$$
consists entirely of BN points.
\end{em}\end{rem}

\bigskip

The translates
of the regions $BGN$ and $M$ can be described explicitly
as follows:
$$ T_{d',s}(BGN)=\{\pa: d'<\mu\leq d'+1,\
0<\lambda\leq\frac{s}{g}(\mu-d'-1)+s,$$
$$\pa\ne(d'+1,s)\}$$
$$    T_{d',s}(M)=\{\pa: d'+1<\mu<d'+2,\ 
0<\lambda\leq\frac{s}{g}(\mu-d'-1)+s\}$$
$$\quad\quad\cup\{(d'+2,\lambda): 0<\lambda<s\}$$
Furthermore, it
follows from these formulae that $T_{d'+1,s}(BGN)$ is
strictly included in  $T_{d',s}(M)$, while
 $T_{d'+1,s+1}(BGN)$ strictly includes  $T_{d',s}(M)$,
i.e. 
$$T_{d'+1,s}(BGN) \subsetneqq T_{d',s}(M) \subsetneqq  
 T_{d'+1,s+1}(BGN).\eqno(5)$$
Therefore, if $d'\geq\hparti$, the translate
$T_{d',s}(BGN)$ covers a larger region than 
does $T_{d'-1,s-1}(M)$.

\bigskip

We will use these relations to translate, in a convenient way, 
the known regions. From Example 3.12, we can assume $g>3$.

\bigskip

We can {\it obtain a new}  region
of  BN points, {\it either by translating} $BGN$ by $T_{2,1}$,
or $M$ by $T_{1,1}$; by $(5)$, the latter
covers a larger area than the former, so we  use
$T_{1,1}(M)$ to enlarge the region.  

 We now {\it continue the process}, translating $M$ by $T_{d',1}$
 with increasing $d'$ (but always keeping $d'<g-2$ to
remain in $R$) (as
 illustrated in Figure 4; there, 
$T_{2,1}(BGN)$ is represented by the lighter
 part in the first diagram, so we can compare it to
$T_{1,1}(M)$).

Now, for exactly the same reasons as before,
this is the best
 we can do {\it as long as} $d'<\hat\eta(2)-1$. However, when
 $d'=\hat\eta(2)-1$, $T_{\hat\eta(2)-1,1}(M)$ covers a
smaller region than $T_{\hat\eta(2),2}(BGN)$, so we now
use the latter (see Figure 5). Of course we can now only
guarantee to get $n$-BN points, for some values of $n$;
however, see Remarks 4.3 and 4.4.

 If we have not yet arrived at $d'=g-2$, at
the {\it next step} we  cover a larger region using now
$T_{\hat\eta(2),2}(M)$. We then {\it continue the process} until
we arrive at $\hat\eta(3)-1$ or $g-2$. In the latter case
{\it we stop}, in the former {\it we use now} $T_{\hat\eta(3),3}(BGN)$,
{\it and repeat the process}.

 The union of trapezia obtained in this way is
therefore a polygonal region which consists entirely of
$n$-BN points for some values of $n$. This is the best 
we can do purely by
translating, but there is a possibility of
{\it obtaining further} $n$-BN points by first translating beyond
$\mu=g-1$ and then {\it applying the Serre duality} map
$\sigma$. 

Thus we consider the affine maps $U_{d',s}=\sigma\circ
T_{d',s}$, given explicitly by
$$U_{d',s}\pa=(2g-2-\mu-d',\ s\lambda+g-1-\mu-d').$$
We now have $d'\ge g-2$ and $T_{d',s}$ maps part of
$BGN\cup M$ below the Riemann-Roch line and hence outside
$P$. So $U_{d',s}\pa$ will not lie entirely in $R$; in
fact the second coordinate in the above formula can be
$\le0$. However it is easy to see that
$$ U_{d',s}(BGN)\cap R\ne\emptyset\iff
d'\ge\hparti\text{ and }g-1\le d'\le\min\{s+g-2,2g-3\},$$
$$ U_{d',s}(M)\cap R\ne\emptyset\iff
d'\ge\hparti\text{ and }g-2\le
d'\le s+g-3.$$
If $s\ge g$, then $\parti\ge s+g-2$, so there are no $d'$
satisfying the above conditions; we shall therefore
assume that $s<g$.

We write for convenience $U_{d',s}(1,1)=(d_1,s_1)$, so
that$$d_1=2g-3-d',\ s_1=s+g-2-d'.$$
Then $U_{d',s}(BGN)\cap R$ is given by
$$d_1\le\mu<d_1+1,\quad0<\lambda\le(1-\frac
sg)(\mu-d_1)+s_1$$with the point $(d_1,s_1)$ omitted,
while $U_{d',s}(M)\cap R$ is given by
$$d_1-1<\mu<d_1,\quad0<\lambda\le(1-\frac
sg)(\mu-d_1)+s_1$$together with the line segment
$$\ell =\{(d_1-1,\lambda): 0<\lambda<s_1-1\}.$$ 

\bigskip

The following lemmas will show that we can gain
 an extra triangle by replacing
$T_{\hat\eta(s_1+1)-2,s_1}(M)$ by  the appropriate 
$(U_{d',s}(BGN)\cap R)\cup \ell ''$ where $\ell ''$
is a line segment.

\begin{lemma}Suppose $d'\ge\hparti$ and $g-2\le
d'\le s+g-3$. Then $s_1\ge1$, $d_1-1\ge\hat\eta(s_1)$ and
$$U_{d',s}(M)\cap R\subset T_{d_1-1,s_1}(BGN)\cup \ell .$$ 
\end{lemma}

{\it Proof:} Note first that
$s_1=s+g-2-d'\ge1$. 

The inequality $d_1-1\ge\hat\eta(s_1)$
is equivalent to $$s_1(d_1-1)\ge(s_1-1)(s_1+g),$$ or,
substituting for $s_1$, $d$,
$$(s+g-2-d')(2g-4-d')\ge(s+g-3-d')(s+2g-2-d').$$
This simplifies to
$$(s+1)d'\ge s^2+(g-1)s-2=(s-1)(s+g)+g-2.$$
But by hypothesis $sd'\ge(s-1)(s+g)$ and $d'\ge g-2$.
This proves the inequality.

Comparing the formulae for
$U_{d',s}(M)$ and $T_{d_1-1,s_1}(BGN)$, we see that it is
now sufficient to prove that $\frac{s_1}g\le1-\frac sg$,
i.e.~$s_1+s\le g$. Now $s_1+s=2s+g-2-d'$, so we need to
show that $d'\ge2s-2$. Since $s<g$, this follows from the
hypothesis $d'\ge\hparti$.\vskip 10pt

Note that, if $s_1=1$, then $\ell =\emptyset$, while, if
$s_1>1$,
$$\ell \subset T_{d_1-2,s_1-1}(BGN).$$
Combined with the lemma, this tells us that $U_{d',s}(M)$
gives nothing new.

\begin{lemma}
Suppose $d'\ge\hparti$ and $g-1\le
d'\le s+g-2$. If $d_1\ge \hat\eta(s_1+1)$, then
$$U_{d',s}(BGN)\cap R\subset T_{d_1,s_1+1}(BGN)\cup
\ell '$$ where $$\ell '=\{(d_1,\lambda): 0<\lambda<s_1\}.$$  
\end{lemma}

{\it Proof:} This follows easily from the formulae for
the two sets.
\vskip 10pt If $s_1=0$, then $\ell '=\emptyset$, while, if
$s_1\ge1$,$$\ell '\subset T_{d_1-1,s_1}(BGN).$$So again we
get nothing new.

It remains therefore to consider the case
where $d_1<\eta(s_1+1)$. 
 Now
$(d_1+1,s_1+1)=\sigma(d',s)$. Since the Brill-Noether
number $\rho$ is invariant under $\sigma$ and
$d'\ge\parti$, it follows that $d_1+1\ge\eta(s_1+1)$. So
the only case we need to consider is
$$d_1+1=\hat\eta(s_1+1)\le g-1.$$

\begin{lemma}In the above circumstances,
$$T_{d_1-1,s_1}(M)\subset (U_{d',s}(BGN)\cap R)\cup
\ell ''$$ where $$\ell ''=\{(d_1+1,\lambda):
0<\lambda<s_1\}.$$
\end{lemma}

{\it Proof:} From the formulae for the two sets, we see
that it is sufficient to prove that  $\frac{s_1}g\le1-\frac
sg$. For this, see the proof of Lemma 4.5.

\bigskip

Thus in each chain
$$T_{\hat\eta(s_1),s_1}(M), \ldots, 
T_{\hat\eta(s_1+1)-2,s_1}(M)$$ in the construction described earlier, we
can gain an extra triangle by replacing\linebreak
$T_{\hat\eta(s_1+1)-2,s_1}(M)$ by  the appropriate 
$(U_{d',s}(BGN)\cap R)\cup \ell ''$.

\bigskip

Finally, then, we {\it define} $BMNO$ to
be the union of the trapezia constructed above together
with their Serre duals. The region $BMNO\cap R$ is bounded
from below by
$\lambda=0$, on the sides by $\mu=0$ and $\mu=g-1$, and
from above by the graph of a seesaw-like function $f_g$
defined on the interval $(0,g-1]$ by 
$$f_g(\mu)=\left\{\begin{array}{lll}
\frac{s}{g}(\mu-\lceil\mu\rceil)+s
&\mu\in(\hat\eta(s),\hat\eta(s)+1]\\
\frac{s}{g}(\mu-\lceil\mu\rceil+1)+s
&\mu\in(\hat\eta(s)+1,\hat\eta(s+1)-1]\\
\frac{\hat\eta(s+1)-s}g(\mu-\lceil\mu\rceil+1)+s
&\mu\in(\hat\eta(s+1)-1,\hat\eta(s+1)].
\end{array}\right.$$

We extend $f_g$ to the whole interval $(0,2g-2)$ by
insisting that its graph is invariant under $\sigma$; the
graph of $f_g$ is then the top boundary of $BMNO$.

Figure 6 shows a typical $BMNO$ region. 

We stress the fact that we have to exclude from $BMNO$
 those  points  corresponding to translates of those parts
of the boundaries  of $BGN$ or $M$ which are not included
in the original regions, and we can summarize  as follows:

{\it If $\pa$ lies in or on the polygon defined above,  $\pa$ is
$n$-BN, for many values of $n$, except for $\mu=0$ 
and $\pa=(1,1)$, and possibly for
$\mu\in{\Bbb N}-\{0,1\}$,  $\mu\in(\hparti,\hat\eta(s+1)]$,
$\lambda\ge s$.}

\begin{rem}\begin{em} When $X$ is not hyperelliptic,
 Mercat has proved
recently \cite{M3} that the results of \cite{M1}
extend to the case $\mu=2$. The constructions of
\cite{M3} are the same as those of \cite{M1}, so the
proofs of section 2 still work, except that in
Lemmas 3.6 and 3.7, we should replace the condition
$d'\le2g$ by $d'\le2g-1$. The effect of this is that those
of the points excluded from $BMNO$ as above which arise as
translates of the right-hand boundary of $M$ can be
restored. However those points arising from the left-hand
boundary of $BGN$ cannot be restored. Thus the only points
of the boundary which must be excluded are those of the
form $(\hparti,\lambda)$ with $\lambda>(s-1)(1+\frac
1g)$ and the points $(\hparti+1,s)$ which arise as
translates of $(1,1)$.
\end{em}\end{rem}

\begin{rem}\begin{em} For semistable bundles, the 
results of \cite{BGN} and
\cite{M3} allow us to include both left-hand and
right-hand boundaries of $BGN\cup M$ (and indeed the point
$(1,1)$). So in this case the whole boundary of $BMNO$ can
be included. Moreover one can include the whole of the line
segments $\{(\hparti,\lambda): 0<\lambda\le s\}$.
\end{em}\end{rem}

\begin{rem}\begin{em} The analysis in Remarks 4.3 and 4.4 
 works also for $U_{d',s}$ and hence
whenever $\lambda<f_g(\mu)$ (with the usual
exceptions for integral values of $\mu$). So there is 
certainly a dense subset of $BMNO$ consisting of BN points.
\end{em}\end{rem}

\bigskip

  Finally, we have the following 
proposition, showing that
the region $BMNO$ always ``stays close'' to the BN curve: 

\begin{propn} Let 
$$\rho_g(\mu)=\frac{\sqrt{(\mu-g+1)^2+4(g-1)}+\mu-g+1}2$$
denote the function whose graph is the BN curve. Then, for
$\mu\in(0,2g-2)$, $$0\le\rho_g(\mu)-f_g(\mu)<1.$$
\end{propn}

{\it Proof:} Since the graphs of $\rho_g$ and $f_g$ are
both invariant under $\sigma$, it is sufficient to
prove this for $\mu\le g-1$.

For the first inequality we
need to prove that every point of
$BMNO$ lies on or below the BN curve. Since this is
certainly true for points of
$BGN\cup M$ and $\tilde\rho$ is invariant under $\sigma$,
it is sufficient to prove that, whenever $\pa\in BGN\cup
M$ and $d'\ge\parti$,
$$I=\frac1{\lambda}\left(\tilde\rho(T_{d',s}\pa)-
\tilde\rho\pa\right)\ge0.$$

A simple calculation shows that
$$
\begin{array}{lll}
I&= sd'-(s-1)((s+1)\lambda-\mu+g-1)\\
&\ge(s-1)(s-(s+1)\lambda+\mu+1)\\
&\ge(s-1)\left(\mu+\dfrac{s+1}g(1-\mu)\right)\ge0.
\end{array}
$$

For the second inequality, note first that both $\rho_g$
and its derivative $\rho_g'$ are strictly increasing
(this is easy to see either geometrically or by
calculus). It follows from the formulae for $f_g(\mu)$
and the fact that, by definition of $\hparti$,
$$\rho_g(\mu)<s+1\text{ for
}\mu\le\hat\eta(s+1)-1,$$ that it is sufficient to prove
the inequalities
$$\rho_g(\hat\eta(s+1))-\left(
\frac{\hat\eta(s+1)-s}g+s\right)<1$$
$$ \rho_g(\hat\eta(s+1))-\left(s+1-
\frac{s+1}g\right)<1$$
for $s\ge1$ and $\hat\eta(s+1)\le g-1$.

Since $\rho_g'(g-1)=\frac12$ and $\rho_g'$ is strictly
increasing,
$$\rho_g(\hat\eta(s+1))<\rho_g(\hat\eta(s+1)-1)+
\frac12<s+\frac32.$$
On the other hand
$$
\begin{array}{lll}
\dfrac{\hat\eta(s+1)-s}g+s&\ge&
\dfrac{\eta(s+1)-s}g+s\\
&=&\dfrac s{s+1}+s\ge s+\dfrac12,
\end{array}
$$
proving the first of the required inequalities. Also
$\eta(s+1)\le g-1$ implies that $(s+1)^2\le g$; hence
$\frac{s+1}g\le\frac1{s+1}\le\frac12$. Thus
$$s+1-\frac{s+1}g\ge s+\frac12$$and we are done.

\begin{rem}\begin{em} A careful analysis
of this proof shows that the worst cases for
$\rho_g(\mu)-f_g(\mu)$ are as $\mu\to\hat\eta(s+1)-1$
from above. Thus in fact
$$\rho_g(\mu)-f_g(\mu)<\max[\rho_g(\hat\eta(s+1)-1)-s]$$
taken over values of $s\ge1$ for which $\hat\eta(s+1)\le
g-1$, and this inequality is best possible. The best
possible inequality which is independent of $g$ is the
one stated in the proposition.
\end{em}\end{rem}

 Examples of stable bundles which are
outside the range to which the constructions of this section 
apply are given in \cite {BF} and \cite{Mu}, and some
different examples in \cite{BMN}.

\section{ Comparison with Teixidor's region}

We now compare $BMNO$ with the corresponding region $T$
constructed by the results of Teixidor \cite{Te1} and Mercat
\cite{M2}, mainly by means of some examples. 

\smallskip

In the stable case, Teixidor's original
result excluded from $\bn$ the vertical segments of length
1, with upper end at a point
 {\it on\/} the BN curve $\tilde\rho=0$ with
integer coordinates. However Mercat in \cite{M2} removed
this restriction
 except for the topmost point of each segment, although
he needs also to exclude all the points described in the
last sentence of the following theorem, while Teixidor excluded
only those segments whose topmost point lies on the BN
curve. We will quote the results of both as follows:

\begin{thm}{\bf (Teixidor/Mercat)} A point $\pa$
determines a non-empty locus $\wbn$ if any of the
following three conditions holds: 
\item $\tilde\rho(\lceil\mu\rceil,
\lceil\lambda\rceil)\geq0\text{ and }
0\neq\lambda-\lfloor\lambda\rfloor\leq\mu-\lfloor\mu\rfloor$
\item $ 
\tilde\rho(\lfloor\mu\rfloor, \lceil\lambda\rceil)\geq0 
\text{ and }
\lambda-\lfloor\lambda\rfloor>\mu-\lfloor\mu\rfloor$
\item $
\tilde\rho(\lfloor\mu\rfloor, \lfloor\lambda\rfloor)\geq0
\text{ and }\lambda=\lfloor\lambda\rfloor.$

Moreover, under the same conditions, $\bn$ is non-empty
except possibly for points $\pa$ with $\mu$, $\lambda$
integers and $\tilde\rho(\mu-1,\lambda)<0$. 
\end{thm}

\begin{rem}\begin{em} In the semistable case, 
this theorem  is a mere
translation of a result of Teixidor
 (\cite{Te1}, Theorem 1, p. 386) to the $\pa$ language;
note that Teixidor's result is stated for $X$ generic,
but for semistable bundles this automatically implies
the result for any $X$. Observe that 
conditions (1) and (2) in fact define triangles in the
$\pa$-plane, with all their vertices at points with integer
coordinates, as illustrated in Figure 7, where
the lighter area corresponds to the first condition and
the darker to the second. Condition 3 describes a 
horizontal
segment of length 1, starting at  the point
$(\lfloor\mu\rfloor, \lfloor\lambda\rfloor)$. 
\end{em}\end{rem}

As shown in Figure 7, for any point with integer coordinates
on or below
the  BN curve, the first two conditions together determine
a parallelogram; hence,
 the region defined by Theorem 5.1 is sometimes
referred to as  ``Teixidor's parallelograms''. We denote 
this region  by $T$.

   We
can describe the region $T$ in a similar way to $BMNO$ by first
defining, for any integer $s$,
$$\hparti'=\lceil\parti+\tfrac1s\rceil-1.$$
Then $d'\ge\hparti'$ if and only if $\tilde\rho(d'+1,s)
\ge0$. (Recall that $d'\ge\hparti$ if and only if
$\tilde\rho(d',s)\ge-1$.)
The region $T$ is then bounded below by $\lambda=0$, on
the sides by $\mu=0$ and $\mu=2g-2$ and from above by the
graph of a function $t_g$ defined by
 $$t_g(\mu)=\cases
\mu-\lceil\mu\rceil+s &\mu\in
(\hparti',\hparti'+1]\\
s&\mu\in(\hparti'+1,\hat\eta(s+1)'] \endcases$$

  Unlike $f_g$, the function $t_g$ is in fact
 continuous and non-decreasing, so the shape of $T$ is
simpler than that of $BMNO$.
Note also that the region covered by Teixidor's
parallelograms is invariant under $\sigma$, so we do not
obtain anything new by using Serre duality. Finally it
is easy to check that $0\le\rho_g(\mu)-t_g(\mu)<1$
(compare Proposition 4.11). 

Figure 8 shows a typical Teixidor polygon
(here, $g=10 $ and the only vertex on the BN  
curve is $(3,9)$, since
$3$ is the only divisor of $g-1=9$).

\bigskip

To compare the
 upper boundaries of $T$ and $BMNO$, we first note that
$$\hparti'=\cases\hparti\text{ if
}\hparti=\parti\\\hparti-1\text{ otherwise}.\endcases$$
For $\mu\le g-1$, it follows that $f_g(\mu)\ge t_g(\mu)$
except possibly in the intervals $(\hparti-1,\hparti+1)$.
If $\hparti=\parti$ (or equivalently
$\tilde\rho(\hparti,s)=-1$), then $f_g(\mu)\ge t_g(\mu)$
in this interval as well. On the other hand, if
$\hparti\ne\parti$, then $t_g(\mu)>f_g(\mu)$ on 
$(\hparti-1,\hparti+1)$. Thus $BMNO$ always extends
outside $T$ and, for almost all values of $g$, $T$ also
extends outside $BMNO$.

 At any rate, for a given (small) genus, it is
easy to  compute both $\hparti$ and $\hparti'$ explicitly.
The figures 9, 10 and 11, illustrate the
 cases $g=10$, $g=12$, and $g=13$, respectively, where different
situations can be appreciated.
 There the shaded area is $BMNO$, and Teixidor's polygons
are only outlined.

\section{ The hyperelliptic case}

Suppose now that $X$ is a non-singular hyperelliptic
curve of genus $g\geq 3$. If we denote by $L$ the hyperelliptic
line bundle on $X$ then $h^0(L^{\otimes(s-1)})=s$ for
$1\le s\le g$, so we can take $d'=2s-2$ in Theorem 3.9.
The analogue of Theorem 4.2 is

\begin{thm} Let $X$ be a non-singular hyperelliptic
curve of genus $g\geq 3$. If $\pa\in
BGN\cup M$ and $1\le s\le g$, then $T_{2s-2,s}\pa$ is
$n$-BN for all $n$ 
such that $(\mu ,\lambda )$ is $n$-BN.
\end{thm}

We now define
 $$BMNO_h=\bigcup_{1\le s\le
g-1}\left(T_{2s-2,s}(BGN\cup M)\cap P\right).$$
It will be convenient to include the point $(2,1)$
in $M$ (see \cite{M2}).

This region is already invariant under Serre duality, so
we do not need to invoke the transformations $U_{d',s}$
in this case. The top boundary of $BMNO_h$ is given by
the graph of the function $h_g$ defined on $(0, 2g-2)$ by
$$h_g(\mu)=\frac{s}{g}(\mu-2s+1)+s\text{ for
}\mu\in(2s-2,2s].$$

The analogues of Remarks 4.3 and 4.4 hold
and indeed we can improve Remark 4.4 (ii). For $1\le s\le
g-1$, the region
$$2s-1<\mu\le2s,\ \ \ \lambda\le s$$
consists entirely of BN points. By Serre duality, so also
does
$$2g-2-2s\le\mu<2g-1-2s,\ \ \ \lambda\le s+\mu-g+1,$$ i.e.
(replacing $s$ by $g-s$)
$$2s-2\le\mu<2s-1,\ \ \ \lambda\le \mu-s+1.$$
Of course, all points of $BGN\cup M$ are BN, hence also
all points of its Serre dual. These results are 
illustrated in Figure 12.

 In the semistable
case, we can include the points $(2s-1,s)$ and also
the line segments $\{(2s,\lambda):s<\lambda\le s+1\}$.

\bigskip

The next step is to show that all special stable
bundles, except for certain line bundles, lie in
$BMNO_h$. 

\begin{thm} Let 
$X$ be a hyperelliptic curve, $E$ a stable bundle on $X$
of rank $n$, degree $d$ and slope $\mu=\frac{d}{n}$, and
$s$ an integer. \item {\it 1)} If $0\le s\le g$
and $2s-2<\mu<2s$, then
$$h^0(E)\le sn+\frac{s}{g}(d-(2s-1)n).$$
\item {\it 2)}If $0\le s\le g-1$, $\mu=2s$ and $E\not\cong
L^{\otimes s}$, then
$h^0(E)\le sn$.
\end{thm}

{\it Proof:} (1) We begin by writing
$$F_s(n,d)=sn+\frac{s}{g}(d-(2s-1)n).$$
We check easily that
$$2F_s(n,d)=F_{s-1}(n,d-2n)+F_{s+1}(n,d+2n).$$

To prove the theorem, we argue by induction on $s$. For
$s=0$, the result is obvious, since $E$ stable with
$\mu<0$ implies $h^0(E)=0$. The result for $s=g$ follows
from this by Serre duality and Riemann-Roch.

Now suppose $0<s<g$. Suppose that there exists a stable
bundle $E$ of slope $\mu$ with $2s-2<\mu<2s$ and such
that $H^0(E)=F_s(n,d)+b_0$ with $b_0>0$. Tensoring the
exact sequence $0\to L^*\to H^0(L)\otimes{\cal O}\to L\to
0$ by $E$, we get
$$0\to L^*\otimes E\to H^0(L)\otimes E\to L\otimes
E\to0.$$ Since $h^0(L)=2$, this gives
$$2h^0(E)\le h^0(E\otimes L^*)+ h^0(E\otimes L).$$
By inductive hypothesis, we have
$$h^0(E\otimes L^*)\le F_{s-1}(n,d-2n);$$
hence
$$h^0(E\otimes
L)\ge2F_s(n,d)+2b_0-F_{s-1}(n,d-2n)=F_{s+1}(n,d+2n)+2b_0.$$
Thus $h^0(E\otimes L)=F_{s+1}(n,d+2n)+b_1$, with
$b_1\ge2b_0$. Continuing in this way, we construct a
sequence $(b_i)$, defined by
$$h^0(E\otimes L^{\otimes i})=F_{s+i}(n,d+2in)+b_i,$$
with
$$b_{i+1}\ge 2b_i-b_{i-1}.$$
We deduce that this sequence is strictly increasing.

On the other hand, by the result for $s=g$, we have
$$h^0(E\otimes L^{\otimes (g-s)})\le F_g(n,d+2(g-s)n).$$
So $b_{g-s}=0$, which is a contradiction. The result
follows.

\bigskip
(2) \ \  Again we proceed by induction. For $s=0$, the only
stable bundle of slope $0$ with $h^0(E)>0$ is $\cal O$.
Similarly, the only stable bundle of slope $2g-2$ with
$h^0(E)>(g-1)n$ is $K$.

For $0<s<g-1$, we proceed as in (1). If there exists a
stable bundle $E$ of slope $2s$ such that $h^0(E)=sn+b_0$
with $b_0>0$, we define the sequence $(b_i)$ for $1\le
i\le g-s-1$ by $h^0(E\otimes L^{\otimes i}) = (s+i)n+b_i$
and prove that $(b_i)$ is strictly increasing. On the
other hand, since by hypothesis $E\otimes L^{\otimes
(g-s-1)}\not\cong K$, it follows that $b_{g-s-1}=0$.
Again we have a contradiction.

\begin{rem}\begin{em} It follows from the
proof of Theorem 6.2 that, if $1\le s\le g-1$ and 
$h^0(E)$ takes its maximum value $F_s(n,d)$ (or $sn$),
then also $h^0(E\otimes L^*)=F_{s-1}(n,d-2n)$ (or
$(s-1)n$) and $h^0(E\otimes L)=F_{s+1}(n,d+2n)$ (or
$(s+1)n$).
\end{em}\end{rem}

\begin{cor}If $\pa \in BMNO_h$, then $\pa$ is $n$-BN 
 for infinitely many values of $n$. The only 
special stable bundles
which lie outside $BMNO_h$ are the line bundles 
$L^{\otimes(s-1)}$
 for $1\leq s\leq g$ and $L^{\otimes(s-1)}(p)$ 
for $1\leq s\leq g-1$
 and $p \in X.$
\end{cor}

{\it Proof:} By Theorems 6.1 and 6.2, it is 
sufficient to prove that
the points $(2s-1,s)$ are only $1$-BN. By 
\cite{BGN} Theorem B, $(1,1)$ is only $1$-BN;
 hence, by Remark 6.3, $(2s-1,s)$ is also 
only $1$-BN.

\bigskip

According to this Corollary, there do not exist 
stable bundles of rank $n>1$ and slope $2s-1$ 
with $1\leq s\leq g-1$ and $h^0(E)=sn.$ However

\begin{propn} Let $X$ be a hyperelliptic
curve. For any integers $n$, $s$ with $n>0$, $1\le s\le
g-1$, there exist stable bundles $E$ of rank $n$ and
slope $2s-1$ with $h^0(E)=sn-1$.
\end{propn}

{\it Proof:} For $s=1$, this is a special case of
\cite{BGN}, Theorem B. If $1<s\le g-1$, a result of
\cite{BMN} says that, if $\Delta$ is a torsion sheaf of
length $n$ with support $n$ distinct points of $X$, and
if $M$ is a line bundle of degree $2$ on $X$ such that
$h^0(M)=1$ then a sufficiently general extension
$$0\to L^{\otimes(s-1)}\oplus\cdots\oplus
L^{\otimes(s-1)}\oplus L^{\otimes(s-2)}\otimes M\to
E\to\Delta\to0$$ is stable, and clearly $h^0(E)=sn-1$.

\bigskip

 We have now completely settled the
nonemptiness problem for bundles of integral
slope. For bundles of non-integral slope, however, we
still have an indeterminate region of points which we
know to be $n$-BN but which may fail to be BN. The next 
example shows that this can indeed
happen.

\begin{exa}\begin{em} Suppose that 
$X$ has genus
$g\geq 4$. Suppose that
$1\le s\le g-1$ and that $E$ is a stable bundle of
rank $n$ and degree $d$ with
$2s-1<\mu=\frac{d}{n}<2s$. Write
$$d-(2s-1)n=gl+l'\text{ with }0\le l'<g.$$
By Theorem 6.2, we have
$$h^0(E)\le
sn+\frac{s}{g}(d-(2s-1)n)=sn+sl+\frac{sl'}{g},$$
in other words
$$h^0(E)\le sn+sl+\lfloor\frac{sl'}{g}\rfloor.$$
If $\lfloor\frac{sl'}{g}\rfloor<1$, then Theorem 6.1
gives the existence of a bundle $E$ with the maximum
possible number of sections.

Suppose now that $s=2$ and $\frac{g}2\le l'<g-1$. 
We claim that, in this case,
$$2h^0(E)<2n+2l+\lfloor\frac{2l'}{g}\rfloor=2n+2l+1.$$

\vskip 6pt{\it Proof of the claim:} Suppose that there exists a stable
bundle $E$ as above with $h^0(E)=2n+2l+1$. We know that
$$2h^0(E)\le h^0(E\otimes L)+h^0(E\otimes L^*).$$
Hence
$$4n+4l+2\le3n+3l+\lfloor\frac{3l'}{g}\rfloor+n+l.$$
So $\lfloor\frac{3l'}{g}\rfloor=2$ and $h^0(E\otimes
L)=3n+3l+2$. Beginning again with $E\otimes L$ and
continuing in this way for a total of $g-3$ steps, we
obtain
$$\lfloor\frac{(g-1)l'}{g}\rfloor=g-2,$$
hence $l'=g-1$. This contradicts our assumption and
proves that there are points which fail to be BN.
\end{em}\end{exa}

\begin{rem}\begin{em} In the exceptional
case $l'=g-1$ of Example 6.6, we can prove that $E$ does
exist. In fact, since $3<\mu<4$, by \cite{BGN} we can
find a stable bundle $F$ of rank $n$ and slope $4-\mu$
with $h^0(F)=n-l-1$. Then $K\otimes F^*$ has slope
$2g-6+\mu$ and
$$h^0(K\otimes
F^*)=(2g-3)n+gl+l'-(g-1)n+n-l-1=(g-1)n+(g-1)l+g-2.$$
Now take $E=K\otimes F^*\otimes
L^{*\otimes(g-3)}=F^*\otimes L^{\otimes2}$ and use the
argument of Example 6.6 in reverse. We obtain
$h^0(E)=2n+2l+1$ as required.
\end{em}\end{rem}

Re's improvement of the Clifford bound for $X$
non-hyperelliptic \cite{R} is intriguingly close to the
boundary of $BMNO_h$. The results of this section show
the extent to which Re's bound fails for a hyperelliptic
curve.

We finally remark that, in the hyperelliptic case, the upper boundary
of the region where $n$-BN points exist is not the graph of a
continuous function; possibly this extends to other cases.

\vspace{1 true cm}
{Fig. 1} The BN map

\vspace{1 true cm}
{Fig. 2} Teixidor's parallelograms

\vspace{1 true cm}
{Fig. 3} The regions $BGN$ and
$M$

\vspace{1 true cm}
{Fig. 4} First steps in the
construction of the region $BMNO$

\vspace{1 true cm}
{Fig. 5} Gain by translating $BGN$ with bundle with 2 sections

\vspace{1 true cm}
{Fig. 6} Construction of a typical $BMNO$ region 
(genus 10)

\vspace{1 true cm}
{Fig. 7} Teixidor's triangles
 
\vskip 10pt

\vspace{1 true cm}
{Fig. 8}
Teixidor's region for genus 10

\vspace{1 true cm}
{Fig. 9} $BMNO$ and $T$ regions for genus $10$

\vskip 10pt
 
\vspace{1 true cm}
{Fig. 10}
$BMNO$ and $T$ regions for genus $12$ (restricted to $R$)

\vspace{1 true cm}
{Fig. 11} $BMNO$ and $T$
regions for genus $13$ (restricted to $R$)

\vspace{1 true cm}
{Fig. 12}
The hyperelliptic case.


\end{document}